\documentclass[12pt]{article}

\usepackage{amsmath}
\usepackage{amssymb}
\usepackage{amsfonts}
\usepackage{amsthm}

\newtheorem{theorem}{Theorem}

\newtheorem{lemma}[theorem]{Lemma}

\theoremstyle{definition}
\newtheorem{defn}       [theorem] {Definition}

\theoremstyle{definition}
\newtheorem*{remark}	{Remark}

\newcommand{\eqn}[2]{\begin{equation}\label{#1}#2\end{equation}}
\newcommand{\eqnst}[1]{\begin{equation*}#1\end{equation*}}

\newcommand{\eqnsplst}[1]{\begin{equation*}\begin{split}%
	#1\end{split}\end{equation*}}

\def\deg{\mathrm{deg}}

\def\Zd{{{\mathbb Z}^d}}
\def\Z2{{{\mathbb Z}^2}}
\def\Z{{\mathbb Z}}

\def\Lam{\Lambda}

\def\muL{\mu^{\rm L}}
\def\muR{\mu^{\rm R}}
\def\muS{\mu^{\rm S}}
\def\OmegaL{\Omega^{\rm L}}
\def\OmegaR{\Omega^{\rm R}}
\def\OmegaS{\Omega^{\rm S}}
\def\partialL{\partial^{\rm L}}
\def\partialR{\partial^{\rm R}}
\def\sigmaL{\sigma^{\rm L}}
\def\sigmaR{\sigma^{\rm R}}

\def\eps{\varepsilon}
\def\es{\emptyset}

\def\cA{{\cal A}}

\def\cC{{\cal C}}

\def\cP{{\cal P}}

\def\bOmega{{\overline{\Omega}}}

\def\Cmax{C^{\mathrm{max}}}
\def\fmax{f^{\mathrm{max}}}
\def\htop{h_{\mathrm{top}}}

\def\PTRFsh{Probab.~Theory Related Fields}

\def\CMPsh{Commun. Math. Phys.}

\begin{document}

\title{Ladder Sandpiles}
\author{Antal A.~J\'arai\footnote{Carleton University, 
School of Mathematics and Statistics, 1125 Colonel By Drive,
Ottawa, ON \ K1S 5B6, Canada.
E-mail: {\tt jarai@math.carleton.ca}. Partially supported by
NSERC of Canada.} \ and
Russell Lyons\footnote{Dept.\ of Mathematics, Indiana University,
Bloomington, IN 47405-5701 USA.
E-mail: {\tt rdlyons@indiana.edu}.
Partially supported by NSF Grant DMS-0231224.}}
\date{}
\maketitle

\begin{abstract}
We study Abelian sandpiles on graphs of the form 
$G \times I$, where $G$ is an arbitrary finite 
connected graph, and $I \subset \Z$ is a finite 
interval. We show that for any fixed $G$ with 
at least two vertices, the stationary measures 
$\mu_I = \mu_{G \times I}$ have two extremal
weak limit points as $I \uparrow \Z$. The extremal
limits are the only ergodic measures of maximum 
entropy on the set of infinite recurrent configurations.
We show that under any of the limiting measures,
one can add finitely many grains in such a way that almost surely
all sites topple infinitely often.
We also show that the extremal limiting measures
admit a Markovian coding.
\end{abstract}

\section{Introduction}
\label{sec:intro}

The sandpile model was introduced by Bak, Tang and 
Wiesenfeld \cite{BTW87,BTW88}, who used it to illustrate 
the idea of self-organized criticality \cite{Bak96},
a concept that became influential in theoretical physics
\cite{Jensen00}. 
The name Abelian sandpile model (ASM) was coined by 
Dhar \cite{Dhar90}, who discovered its Abelian 
property. The model also appeared independently in the 
combinatorics literature, where it is known 
as the chip-firing game, introduced in \cite{BLS91}. 
One of the remarkable features of the ASM is that its simple 
local rules give rise to complex long-range dynamics.
See \cite{Dhar06} for an overview.

Recently, a number of papers were devoted to sandpiles on infinite
graphs, obtained as limits of sandpiles on finite subgraphs
\cite{MRSvM00,MRS02,MRS03,AJ04,JR06}. See the reviews
\cite{MRS05,Jarai05,Redig05} for the main ideas of these 
developments.
A natural approach to studying sandpiles in infinite volume 
is the following. 
Start with the stationary measure $\mu_V$ of the model 
on a finite subgraph $V$,
and characterize the set of weak limit points of 
$\{ \mu_V : \text{$V$ finite} \}$. Then study avalanches
on infinite configurations under the limiting measures,
and construct a dynamics for the infinite system, if 
at all possible. The papers mentioned carry out this program 
to various degrees for the following infinite graphs:
$\Z$; an infinite regular tree; $\Zd$, $d \ge 2$ with or 
without dissipation; and finite-width strips with dissipation. 
Already the first step, determining the
limiting behaviour of $\mu_V$, is usually non-trivial. 
For each infinite graph mentioned above, there is a unique limit
point, and a number of different techniques have been 
employed to show this. Currently there seems to be no 
unified method that applies to a general infinite graph. 

In this paper, we study the sandpile measures on graphs 
that are the product of a finite connected graph $G$ 
with a finite interval $I \subset \Z$,
with particular view towards the limiting behaviour as
$|I| \to \infty$. We call these ``ladder graphs", where 
the ``rungs" of the ladder consist of copies of $G$.
The only dissipative sites in our model are the ones at 
the two ends of the ladder, that is, the sites in
$G \times \{ \text{endpoints of $I$} \}$. Hence, even when 
$G$ is a finite interval in $\Z$, our models are different 
from the dissipative strips studied in \cite{MRS03}.

When $G$ is a single vertex, the model lives on an interval 
$I \subset \Z$. In this well-known case, $\mu_I$ can be found 
explicitly, and its limiting behaviour is trivial: the limit is 
concentrated on a single configuration with constant 
height $2$. However, as was observed in \cite{AD95},
the behaviour of this model is atypical of general 
one-dimensional sandpiles. Already the simplest modifications, 
such as the ``decorated chains" studied in \cite{AD95}, 
give rise to non-trivial limits with positive entropy.

Not surprisingly, the limiting behaviour is also non-trivial 
in our case, provided $G$ has more than one vertex. As we show in 
Theorem \ref{thm:therm} in Section \ref{ssec:weak-limits},
when $|G| \ge 2$, the set of weak limit points consists
of all convex combinations of two different extremal 
measures $\muL$ and $\muR$, which are related by a
reflection of $\Z$. These measures arise from restricting the 
burning algorithm to act exclusively from the left or 
the right, and we call them the left- and right-burnable 
measures. In the case of $|G| = 1$, the left- and right-burnable
measures happen to coincide.

In Theorem \ref{thm:max-entropy} in 
Section \ref{ssec:weak-limits} we show that $\muL$ and
$\muR$ are the only two ergodic measures of maximum 
entropy on the set of infinite recurrent configurations,
and that there is a unique measure $\muS$ of maximum
entropy that is invariant under reflection.
It would be interesting to see whether there is a unique 
measure of maximum entropy for the infinite graphs studied
earlier (where the weak limits are unique).

For most of our arguments on weak limits, some quite general
properties of the model are sufficient. For example, existence
of the limit of left- and right-burnable measures follows from 
the existence of renewals: if all sites in $G \times \{i\}$ 
have the maximum possible height for a fixed $i \in I$, then 
the subconfigurations to the left and right of $i$ are
conditionally independent.

Given an infinite configuration, we can ask what happens if particles 
are added and then the configuration is relaxed.
First, in the case of $\Z$, it is easy to see that if we add a single 
grain to the system (having constant height $2$), then every site topples
infinitely often. 
On the ladder $G \times \Z$ with $|G| \ge 2$, 
finite avalanches do occur with positive probability. However, as 
we show in Section \ref{sec:avalanches}, it is possible to 
add a fixed number of grains in such a way that almost surely 
every site topples infinitely often, with respect to 
any of the limiting measures. 
Hence, there is no sensible dynamics for $G \times \Z$ in general.

{\bf Open question.} \emph{Does the probability of infinitely many 
topples at $(0, 0)$, when adding 1 grain to $(0, 0)$, 
tend to $0$ for $G = \Z_n$ and $n \to\infty$? Here $\Z_n$ is 
the cycle of length $n$.}

The measure $\muL$ (and $\muR$) can be regarded as a subshift on 
a finite alphabet (that depends on $G$), by grouping sites on each 
copy of $G$ together. The set of recurrent configurations are 
characterized in terms of finite forbidden words, and there
is an infinite number of constraints. As our results show (see 
Lemma \ref{lem:entropy-compare}), the number of constraints
grows at a rate strictly smaller than the topological entropy.
Hence the set of left-burnable configurations is a subshift of
quasi-finite type, in the terminology of \cite{Buzzi05}.
In fact, our subshift turns out to be more special. As we show in 
Section \ref{sec:coding}, it admits a Markovian coding, and 
hence it is a sofic shift \cite[Theorem 3.2.1]{LM95}.
We note that the set of all recurrent configurations is also a 
sofic shift (by arguments similar to those for
Lemma \ref{lem:entropy-compare}). However, since recurrent 
configurations lead to mixtures of $\muL$ and $\muR$, we 
study only the latter in detail.
An alternative approach to our results in 
Section \ref{sec:recurrent} would be to analyze the
Markovian coding obtained in Section \ref{sec:coding}. 
However, we prefer to present more direct arguments.

We will assume throughout that the reader is familiar with the basic 
properties of the ASM that can be found in 
\cite{Dhar06,Dhar99a,Dhar99b,IP98,MRZ02}.

\section{Models considered}
\label{sec:models}

Throughout, $G$ will be an arbitrary fixed finite connected graph. 
For $n \le m$, let $I_{n,m}$ denote the graph on the vertex set 
$\{n, \dots, m\}$ with nearest neighbour edges. Let 
$\deg_G(x)$ denote the degree of a vertex $x$ in $G$. We consider 
Abelian sandpiles \cite{MRZ02} defined on the product graph 
$\Lam_{n,m} := G \times I_{n,m}$ (whose edges join vertices $(x, k)$ 
and $(y,\ell)$ when either $x \sim y$ in $G$ and $k=\ell$ or $x= y$ 
and $|k-\ell| = 1$). We are primarily interested in the limit sandpiles
as $n \to -\infty$ and $m \to \infty$ that live on the graph
$\Lam := G \times \Z$.
For convenience, we also introduce $\Lam_{-\infty,m}$ and
$\Lam_{n,\infty}$ with the obvious meaning. We refer to 
$G \times \{ k \}$ as the \emph{rung} at $k$. 

We let $\Delta$ denote the graph Laplacian on $\Lambda$, that is,
the following matrix indexed by vertices in $\Lam$:
\eqnst
{ \Delta_{uv} 
  := \begin{cases}
    \deg_G(x) + 2 & \text{if $u = v = (x,k)$;}\\
    -1            & \text{if $u$ and $v$ are neighbours;}\\
    0             & \text{otherwise.}
    \end{cases} }
For finite vertex-subsets $V \subset \Lam$ that induce a connected
subgraph, we let 
$\Delta_V$ denote the restriction of $\Delta$ to the pairs
$(u,v) \in V \times V$.	In other words, a sink site is added to 
$V$, and each $u \in V$ is connected to the sink by 
$\Delta_{uu} - \deg_V(u)$ edges.

We are interested in the sandpile with toppling matrix 
$\Delta_{\Lambda_{n,m}}$. We will study the case when $G = I_{0,1}$
quite explicitly for illustration. 

The space of stable configurations on a set $V \subset \Lam$ 
is 
\eqnst
{ S_V 
  := 
  \prod_{u \in V} \{1, \dots, \Delta_{u u} \}. } 
We write $S := S_\Lam$. 
For a convenient notation, we define $m(x) := \deg_G(x) + 2$, 
$x \in G$, which is the maximum allowed height at a site 
$u = (x,k)$.
We write $\Omega_V$ for the set of recurrent 
configurations \cite{MRZ02} on $V$ when $V$ is finite.
We define
\eqnsplst
{ \Omega
  &:= \Omega_\Lam
  := \{ \text{recurrent configurations on $\Lam$} \} \\
  &:= \{ \eta \in S : \text{$\eta_W \in \Omega_W$ for
     all finite $W \subset \Lam$} \}. }

For $V \subset \Lam$, if $\Lam \setminus V$ has a
connected component fully infinite to the left (that is,
containing $\Lam_{-\infty,m}$ for some $m$), we denote
that connected component $V^-$. We similarly define 
$V^+$ to the right (which may coincide with $V^-$). 
We define the \emph{left (interior) boundary} of $V$ as
\eqnst
{ \partialL_0 V
  := \left\{ v \in V :\ \text{$v$ has a neighbour in 
    $V^-$} \right\}. } 
We define $\partialR_0 V$ analogously.

\section{Description of recurrent configurations}
\label{sec:recurrent}

\subsection{Left- and right-burnable measures}
\label{ssec:left-right-burnable}

We define a one-sided version of the burning 
algorithm \cite{MRZ02}.

\begin{defn}
\label{defn:left-burn}
Let $V \subset \Lambda$ be finite.
A configuration $\eta \in S_V$ is called \emph{left-burnable} 
if there is an enumeration $v_1, \dots, v_{|V|}$ of $V$ such
that 
\begin{itemize}
\item[(i)] $v_i \in \partialL_0 
  (V \setminus \{ v_1, \dots, v_{i-1} \})$, $1 \le i \le |V|$;
\item[(ii)] $\eta(v_i) > 
  \Delta_{v_i v_i} 
  - | \{ u \in (V \setminus \{ v_1, \dots, v_{i-1} \})^- : 
  u \sim v_i \} |$.
\end{itemize}
Note that this is the usual burning rule with the restriction that
only sites in the left boundary can be burnt. When $\Lam \setminus V$
is connected, the rule becomes identical to the usual burning rule.
We denote by $\OmegaL_V$ the set of left-burnable configurations on $V$.
We define right-burnable configurations and $\OmegaR_V$ analogously.
\end{defn}

\begin{lemma}
\label{lem:restrict}
Let $V \subset \Lambda$ be finite.
We have $\OmegaL_V \subset \Omega_V$. 
If $\eta \in \OmegaL_V$ and $W \subset V$, then $\eta_W \in \OmegaL_W$.
The same holds for $\OmegaR_V$.
\end{lemma}

\begin{proof}
The sequence $v_1, \dots, v_{|V|}$ required by 
Definition \ref{defn:left-burn} is a valid burning sequence in the 
ordinary burning algorithm, since
\eqnst
{ |\{ u \in (V \setminus \{ v_1, \dots, v_{i-1} \})^- : 
    u \sim v_i \}| 
  \le |\{ u \in (V \setminus \{ v_1, \dots, v_{i-1} \})^c : 
    u \sim v_i \}|. }
Therefore, $\OmegaL_V \subset \Omega_V$. For $W \subset V$,
let $w_1,\dots,w_{|W|}$ be the enumeration of $W$ in the order
inherited from the enumeration of $V$. Since
\eqnst
{ |\{ u \in (V \setminus \{ v_1, \dots, v_{i-1} \})^- : 
    u \sim v_i \}| 
  \le |\{ u \in (W \setminus \{ v_1, \dots, v_{i-1} \})^- : 
    u \sim v_i \}|, }
this is a valid left-burning sequence for $\eta_W$.
\end{proof}

\begin{defn}
For $V = \Lam$, a configuration $\eta \in S$ is called 
{\it left-burnable\/} if $\eta_W$ is left-burnable for every finite 
$W \subset \Lam$. {\it Right-burnable\/} configurations are defined
analogously.  We write $\OmegaL$ and $\OmegaR$ for
the sets of these configurations.  
We write $\OmegaS := \OmegaL \cap \OmegaR$.
\end{defn}

\begin{defn}
Let $\mu_{n,m}$ denote the uniform measure on the set of 
recurrent configurations on $\Lam_{n,m}$. 
We denote by $\muL_{n,m}$ the uniform measure on 
left-burnable configurations on $\Lam_{n,m}$ and define 
$\muR_{n,m}$ and $\muS_{n,m}$ analogously.
\end{defn}

In order to illustrate some of the results to come,
we explicitly describe left-burnable configurations in the
simplest non-trivial case $G = I_{0,1}$.

\begin{lemma}
\label{lem:left-burn}
Assume $G = I_{0,1}$. A configuration $\eta \in \Omega_{n,m}$
is left-burnable if and only if the following 3 conditions hold:
\begin{enumerate}
\item each rung contains a $3$;
\item if the rung at $k$ is $(3,1)$, then no rung 
  other than $(3,2)$ can occur to the right of $k$ 
  before a $(3,3)$ occurs. That is, the rungs at
  $k, k+1, \dots$ are of the form:
\eqn{e:leftburn}
{ \begin{matrix}
  3 & 3 & \dots & 3 & 3 & \dots \\
  1 & 2 & \dots & 2 & 3 & \dots  
  \end{matrix} }
with the possibility that there is no $(3,2)$ rung
at all, and the exception that the $(3,3)$ may be missing 
if the right end of $\Lam_{n,m}$ was reached;
\item if the rung at $k$ is $(1,3)$, then no rung 
  other than $(2,3)$ can occur to the right of $k$
  before a $(3,3)$ occurs.
\end{enumerate}
The same holds for right-burnable configurations with 
left and right interchanged.
\end{lemma}

\begin{proof}
By symmetry, we may restrict to the left-burnable case.
It is straightforward to verify that a configuration 
satisfying 1--3 in the Lemma is left-burnable. Namely, the 
configuration can be burnt rung-by-rung, except
when a $(3,1)$ or a $(1,3)$ is encountered. In the
latter case, observe that the configuration in 
\eqref{e:leftburn} is left-burnable (as well as the 
one obtained by exchanging the rows).

Assume now that we are given a left-burnable configuration,
and we show that 1--3 hold.
The proof is by induction on the number $N = m - n + 1$ of rungs.
The case $N = 1$ is trivial. Assume now that $N > 1$ and that
the statement holds whenever the number of rungs is
less than $N$. Observe that the leftmost rung has to 
contain a $3$, otherwise the burning cannot start. 

\emph{Case 1. The leftmost rung is $(3,3)$, $(3,2)$
or $(2,3)$.} Then without loss of generality, we may 
assume that the burning starts with removing the leftmost
rung. Since the $N - 1$ remaining rungs are left-burnable,
the induction hypothesis implies the claim.

\emph{Case 2. The leftmost rung is $(3,1)$ or $(1,3)$.}
We may assume the leftmost rung is $(3,1)$. Then, by the
burning procedure, the next rung is of the form 
$(3,z)$. If $z = 3$, we can use the induction hypothesis
for $N - 2$. The value $z = 1$ leads to a forbidden 
subconfiguration $1\ 1$. If $z = 2$, we can iterate 
the present argument until a rung of the form $(3,3)$
is reached, noting that configurations of the form
$1\ 2\ \dots\ 2\ 1$ are forbidden.
\end{proof}

It follows from the description in Lemma \ref{lem:left-burn}
that $(3,3)$ rungs are renewals, that is, given that 
rung $k$ is $(3,3)$, the subconfigurations to the left and
right are conditionally independent for the appropriate measure $\mu_{n,
m}^\bullet$. The analogous statement
holds for maximal rungs on a general graph, and we prove 
this next. 

The following terminology will be useful.
Let $\cC := \cC(G) := \OmegaL_{\Lam_{0,0}}$ denote the 
set of left-burnable configurations on a single rung. We 
claim this is the same as the set of all recurrent configurations 
on $G \times \{ 0 \}$ with at least one $x \in G$ such that 
$\eta(x,0) = m(x)$. Indeed, no burning will occur without such an $x$, and
with such an $x$, we can left-burn $(x,0)$, and once 
this is done, left-burning becomes equivalent to ordinary burning (i.e.,
burning from both sides) since the left and right boundaries merge.
By the same reasoning, $\cC$ is also the set of right-burnable 
configurations at a single rung. 
By abuse of notation, we regard $\cC$ as a set of configurations on any
particular rung.
For $\eta \in \OmegaL_{n,m}$, let
$C_k := C_k(\eta) := \eta_{\Lam_{k,k}}$ denote the rung at $k$, 
which is in $\cC$ by Lemma \ref{lem:restrict}.
Let $\Cmax := \Cmax(G) \in \cC$ denote the configuration on $G$ defined
by $\Cmax(x) := m(x)$, $x \in G$. The configuration $\Cmax$ is the maximal 
configuration that can occur on a rung.

\begin{lemma}
\label{lem:renewals-finite} {\bf (Renewals)}
For the measures $\muL_{n,m}$, $\muR_{n,m}$ and $\muS_{n,m}$,
maximal rungs are renewals, that is, given $C_k = \Cmax$, 
the subconfigurations to the left and right of rung $k$ 
are conditionally independent.
\end{lemma}

\begin{proof}
First consider the left-burnable measure.
Let $\eta \in \OmegaL_{n,m}$, and assume that $C_k(\eta) = \Cmax$.
By Lemma \ref{lem:restrict}, both $\eta_{\Lam_{n,k-1}}$ and 
$\eta_{\Lam_{k+1,m}}$ are left-burnable. We need to show that the
two vary independently, that is, for any $\xi \in \OmegaL_{n,k-1}$ 
and $\zeta \in \OmegaL_{k+1,m}$, we have 
$\eta' = \xi \vee \Cmax \vee \zeta \in \OmegaL_{n,m}$, where $\vee$
indicates concatenation. 
Start left-burning on $\eta'$. Since $\xi$ is left-burnable,
there will be a first time when a site $(x,k-1)$ is
burnt. When this happens, we can fully burn rung $k$.
After rung $k$ is burnt, both the rest of $\xi$ and all of
$\zeta$ can be burnt, because they are left-burnable.
Hence $\eta'$ is left-burnable, and $\muL_{n,m}$ has
the renewal property since it is uniform on $\OmegaL_{n, m}$.

The statement for $\muR_{n,m}$ follows by symmetry. The statement
for $\muS_{n,m}$ can be proved by a very similar argument, now
showing that burning from both left and right can be performed.
\end{proof}

In order to investigate weak convergence of the finite-volume
measures, we are going to use some comparisons between the
growth rates (topological entropies) of certain sets of configurations.
This is formulated in the lemma below. 

Let $a_n := |\OmegaL_{1,n}|$. By Lemma \ref{lem:restrict},
we have 
\eqnst
{ a_{n+m} 
  \le |\OmegaL_{1,n}| |\OmegaL_{n+1,n+m}| 
  = a_n a_m. }
Therefore,
\eqn{e:hL}
{ h^{\rm L} 
  := \lim_{n \to \infty} \frac1n \log |\OmegaL_{1,n}| 
  = \inf_{n \ge 1} \frac1n \log a_n }
exists. The limit $h^{\rm L}$ is the topological entropy $\htop(\OmegaL)$
of $\OmegaL$ with respect to translations \cite{KH95}.
For any $\delta > 0$, there exists $C = C(\delta)$
such that
\eqnst
{ \exp \{ h^{\rm L} n \} 
  \le a_n
  \le C \exp \{ h^{\rm L} n (1 + \delta) \}. }
By symmetry, $h^{\rm L} = h^{\rm R} = \htop(\OmegaR)$.
We also define $s_n := |\OmegaS_{1,n}|$, and again, by
submultiplicativity, we have 
\eqnst
{ h^{\rm S} 
  := \lim_{n \to \infty} \frac1n \log s_n
  = \htop(\OmegaS). }
We further define the spaces
\eqnsplst
{ \Omega^{\rm L,0}
  &:= \{ \eta \in \OmegaL : C_k(\eta) \not= \Cmax,\, 
    -\infty < k < \infty \} \\
  \Omega^{\rm L,0}_{n,m}
  &:= \{ \eta \in \OmegaL_{n,m} : C_k(\eta) \not= \Cmax,\, 
    n \le k \le m \}, }
and we define $\Omega^{S,0}$ and $\Omega^{S,0}_{n,m}$
analogously in the symmetric case.
Let $b_n := |\Omega^{\rm L,0}_{1,n}|$, $r_n := |\Omega^{S,0}_{1,n}|$,
$h^{\rm L,0} := \lim_{n \to \infty} (1/n) \log b_n = \htop(\Omega^{\rm L,0})$ 
and $h^{S,0} := \lim_{n \to \infty} (1/n) \log r_n = \htop(\Omega^{S,0})$.

\begin{lemma}
\label{lem:entropy-compare}
Assume that $G$ is not a single vertex. Then
\begin{itemize}
\item[(i)] $0 < h^{\rm L,0} < h^{\rm L}$; 
\item[(ii)] $0 < h^{S,0} < h^{\rm S}$;
\item[(iii)] $h^{\rm S} < h^{\rm L}$.
\end{itemize}
\end{lemma}

\begin{proof}
For $x \in G$, define the rung 
\eqnsplst
{ C^x (z)
  := \begin{cases}
    m(x) - 1 & \text{$z = x$,} \\
    m(z)     & \text{$z \not= x$.} 
    \end{cases} }
It is straightforward to check that since $G$ consists of more than 
one vertex, $C^x \in \cC$. Now let $x, y \in G$, $x \not= y$.
Any sequence consisting exclusively of rungs $C^x$ and $C^y$ is 
both left- and right-burnable. Hence $\log 2 \le h^{S,0} \le h^{\rm L,0}$.

For $0 \le k \le n$ we select $k$ of the rungs. Consider the function
that changes these $k$ rungs of an $\eta \in \Omega^{\rm L,0}_{1,n}$ to 
$\Cmax$. Any configuration so obtained is in $\OmegaL_{1,n}$ and has 
at most $|\cC|^k$ preimages. Therefore, the number of different new
configurations obtained is at least $b_n / |\cC|^k$. Summing over 
$k$ and all choices of $k$ rungs, we have
\eqnst
{ a_n 
  \ge b_n \sum_{k=0}^n \binom{n}{k} \frac1{|\cC|^k}
  = b_n \left( 1 + \frac1{|\cC|} \right)^n. }
Hence, $h^{\rm L} \ge  h^{\rm L,0} + \log(1 + |\cC|^{-1})$. This proves (i).
By a similar argument, $h^{\rm S} \ge h^{S,0} + \log(1 + |\cC|^{-1})$,
which proves (ii). 

The argument to prove (iii) is also similar: Note that when $G = I_{0,1}$,
Lemma \ref{lem:left-burn} implies that the sequence of rungs 
$(3,3)$, $(2,3)$, $(3,1)$, $(3,3)$ is left-burnable, but not right-burnable.
We can adapt this observation to general $G$. Let 
$x \sim y \in G$, and define the rungs
\eqnst
{ C_1 (z)
  := \begin{cases}
    m(z) - 1   & z \not= y, \\
    m(y)       & z = y,
    \end{cases}
  \qquad\qquad
  C_2 (z)
  := \begin{cases}
    m(z) - 1   & z \not= x,y, \\
    m(x)       & z = x, \\
    1          & z = y.
    \end{cases} }

We claim that the sequence $\xi = \Cmax, C_1, C_2, \Cmax$ is left-burnable,
but not right-burnable. In case of left-burning, $\Cmax$ burns first, then
site $y$ burns in $C_1$, and after that the rest of $C_1$ can be burnt.
Now site $x$ burns in $C_2$, and after this the $\Cmax$ rung to the
right burns. This makes the rest of the sites but $y$ in $C_2$ burnable;
finally, site $y$ in $C_2$ can be burnt. In case of right-burning, the 
$\Cmax$ rung on the right can be burnt. After this, site $x$ in $C_2$
can be burnt. This may make other sites in $C_2$ burnable. However, 
crucially, $y$ in $C_2$ cannot be burnt (since it has a neighbour in $C_1$), and
no site in $C_1$ can be burnt, since burning could only start at $y$,
which is ``blocked" by the $1$ in $C_2$.

Assume now that $\eta \in \OmegaS_{1,4n}$, and subdivide $[1,4n]$
into $n$ intervals of length $4$. Consider the mapping that replaces
the rungs at a fixed set of $k$ of these intervals by $\xi$.
The configurations obtained are in $\OmegaL_{1,4n}$, and since
$\xi$ is not right-burnable, they are not in $\OmegaS_{1,4n}$.
The number of preimages of a given element of $\OmegaL_{1,4n}$ 
is at most $|\cC|^{4k}$. Hence we get
\eqnst
{ a_{4n} 
  \ge s_{4n} \sum_{k=0}^n \binom{n}{k} \frac1{|\cC|^{4k}}
  = s_{4n} \left( 1 + \frac1{|\cC|^4} \right)^n. }
This implies (iii).
\end{proof}

\subsection{Weak limits}
\label{ssec:weak-limits}

\begin{lemma}
\label{lem:weak-lim}
The weak limits 
\eqnst
{ \muL := \lim_{\substack{n \to -\infty \\ m \to \infty}} \muL_{n,m},
  \qquad
  \muR := \lim_{\substack{n \to -\infty \\ m \to \infty}} \muR_{n,m},
  \qquad
  \muS := \lim_{\substack{n \to -\infty \\ m \to \infty}} \muS_{n,m} }
exist. The limit measures
$\muL$, $\muR$ and $\muS$ are concentrated on $\OmegaL$,
$\OmegaR$ and $\OmegaS$ (respectively). 
\end{lemma}

\begin{proof}
We first strengthen \eqref{e:hL} to show that with 
$\lambda = \exp(-h^{\rm L})$, the limit
$\lim_{n \to \infty} \lambda^n a_n$ exists and is positive.
Lemma \ref{lem:renewals-finite} implies the renewal equation:
\eqn{e:ren-eqn}
{ a_n
  = b_n + \sum_{k = 1}^n b_{k-1} a_{n-k}, \qquad n \ge 0, }
where we set $a_0 = 1$, $b_0 = 1$. 
Let
\eqnst
{ F(z)
  := \sum_{n = 0}^\infty a_n z^n 
  \qquad \hbox{ and } \qquad
  G(z)
  := \sum_{n = 0}^\infty b_n z^n }
be the generating functions of $\{ a_n \}_{n \ge 0}$ and 
$\{ b_n \}_{n \ge 0}$.
The radius of convergence of $F$ is $\lambda = \exp(-h^{\rm L})$
and that of $G$ is $\exp(-h^{\rm L, 0}) > \lambda$. 
The relation \eqref{e:ren-eqn} implies
\eqnsplst
{ F(z) 
  = \frac{G(z)}{1 - z G(z)}, \qquad 0 \le z < \lambda. }
Since $G$ is analytic in a disc of radius larger than 
$\lambda$, but $F$ has a singularity on the circle
$|z| = \lambda$, we need to have 
$1 = \lim_{z \nearrow \lambda} z G(z) = \lambda G(\lambda)$.
It follows that $p_n := \lambda^n b_{n-1}$, $n \ge 1$, is a 
probability distribution, and with $c_n := \lambda^n a_{n-1}$,
\eqref{e:ren-eqn} has the probabilistic form
\eqn{e:ren-eqn2}
{ c_{n+1}
  = p_{n+1} + \sum_{k = 1}^n p_k c_{n-k}, \qquad n \ge 0. } 
By the Renewal Theorem \cite[page 330]{F68}, we have
$\lim_{n \to \infty} c_n = (\sum_{k \ge 1} k p_k)^{-1}$.
Hence we have
\eqn{e:expasym}
{ \lim_{n \to \infty} \lambda^n a_n
  = \lim_{n \to\infty} c_n/\lambda
  = \frac{1}{\lambda\sum_{k \ge 1} k p_k}
  = \frac{1}{\lambda^2 [\lambda G(\lambda)]'}
  =: \alpha > 0. } 

We are ready to establish the existence of $\muL$. Fix $k \ge 1$
and an elementary cylinder event depending on the rungs $-k,\dots,k$.
That is, we fix $\eta_0 \in \OmegaL_{-k,k}$, and let $E = E(\eta_0)$ 
denote the event that the subconfiguration in rungs $-k,\dots,k$
equals $\eta_0$. We need to show that
\eqnst
{ \lim_{\substack{n \to -\infty \\ m \to \infty}} 
    \muL_{n,m} ( E )
  =: \muL ( E ) \text{ exists.} }
We first show that for $N > k$ large enough and
$-n,m > N$, the event
\eqnst
{ A(N) 
  = \left\{ \text{$\exists$ renewal in $[-N,-k-1]$ and in
    $[k+1,N]$} \right\} }
occurs with high $\muL_{n, m}$-probability. Indeed, letting $\eta$ denote a 
random variable with law $\muL_{n,m}$ and using 
Lemma \ref{lem:entropy-compare}, we get
\eqnsplst
{ &\muL_{n,m} (\text{no renewal in $[k+1,N]$}) 
  = \muL_{n,m} \left( \eta_{\Lam_{k+1,N}} 
    \in \Omega^{\rm L,0}_{k+1,N} \right) \\
  &\qquad \le \frac{|\OmegaL_{n,k}| |\Omega^{\rm L,0}_{k+1,N}|
    |\OmegaL_{N+1,m}|}{|\OmegaL_{n,m}|}
  \le C e^{-\delta (N-k)} }
for some $\delta > 0$ and $C = C(\delta)$ for all large $N$. This implies
\eqn{e:ANbnd}
{ \muL_{n,m} (A(N)^c) 
  \le 2 C e^{-\delta (N-k)}, \qquad -n,m > N. }
On the event $A(N)$, let  
\eqnsplst
{ \tau
  &:= \text{leftmost renewal in $[k+1,N]$} \\
  \hbox{and }\
  \sigma
  &:= \text{rightmost renewal in $[-N,-k-1]$}. }
We also define
\eqnst
{ u (s,t,E)
  := \left| \left\{ \xi \in \OmegaL_{s+1,t-1} : 
    \parbox{2.2in}{$\xi_{\Lam_{-k,k}} = \eta_0$ and no renewal in
    $[s+1,-k-1] \cup [k+1,t-1]$} \right\} \right| }
for $-N \le s \le -k-1$ and $k+1 \le t \le N$.
Considering the values of $\sigma$ and $\tau$ and
counting configurations, we can write
\eqnsplst
{ \muL_{n,m} \big(E,\, A(N)\big)
  &= \sum_{t=k+1}^N \ \sum_{s = -N}^{-k-1}  
     \muL_{n,m} ( E,\, \tau=t,\, \sigma=s ) \\
  &= \sum_{t=k+1}^N \ \sum_{s = -N}^{-k-1} 
     \frac{ a_{s-n} u(s,t,E) a_{m-t} }{ a_{m-n+1} }. }
Using \eqref{e:expasym}, we have
\eqnst
{ \lim_{\substack{n \to -\infty \\ m \to \infty}} 
    \muL_{n,m} \big(E,\, A(N)\big)
  = \alpha \sum_{t=k+1}^N \ \sum_{s = -N}^{-k-1} 
    \lambda^{t-s+1} u(s,t,E). } 
Letting $N \to \infty$ and applying \eqref{e:ANbnd}, we deduce that
\eqnst
{ \lim_{\substack{n \to -\infty \\ m \to \infty}} 
    \muL_{n,m} (E)
  = \alpha \sum_{t=k+1}^\infty \ \sum_{s = -\infty}^{-k-1} 
    \lambda^{t-s+1} u(s,t,E)
  =: \muL (E). }
The statement for $\muR$ follows by symmetry. In the case
of $\muS$, the proof follows a very similar line.
\end{proof}

\begin{remark}
It is not hard to extend the proof above to show that
$\lim_{V \nearrow \Lam} \muL_V = \muL$ (and similarly 
for $\muR$ and $\muS$).
\end{remark}

\begin{lemma}
\label{lem:ergodic}
Maximal rungs are renewals for the measures
$\muL$, $\muR$ and $\muS$, and these measures are ergodic.
If $G$ is not a single vertex, then $\muL$ and $\muR$ are not 
symmetric under reflection, while $\muS$ is. The measures $\muL$ and 
$\muR$ are reflections of each other.
\end{lemma}

\begin{proof}
The renewal property follows from Lemma \ref{lem:renewals-finite}
by passing to the limit. Ergodicity follows from the existence
of renewals. The configuration $\xi$ given in the
proof of Lemma \ref{lem:entropy-compare} (iii) shows that
$\muL \not= \muR$. 
\end{proof}

\begin{remark}
It follows from general arguments that $\mu^L$ and $\mu^R$ have
maximal entropy. For example, by a counting argument one can show that 
$|\Omega_{1,n}|$ has exponential growth rate $h^L$, and this 
allows one to adapt the argument of \cite[Proposition 1.12 (ii)]{BS94}.
Below we show that there are no other measures of maximal entropy. 
\end{remark}

\begin{theorem}
\label{thm:max-entropy}
The only two ergodic measures of maximum entropy on $\Omega$ are 
$\muL$ and $\muR$. The unique symmetric ergodic measure of 
maximum entropy is $\muS$.
\end{theorem}

\begin{proof}
There exists a measure $\mu$ of maximum entropy on $\Omega$ \cite{KH95}. 
By ergodic decomposition, we may assume that $\mu$ is ergodic. 
We show that in this case either $\mu = \muL$ 
or $\mu = \muR$, which shows that these are the only two 
ergodic measures of maximum entropy.

We first show that $\mu \{ C_0 = \Cmax \} > 0$.
To see this, note that increasing the height of any site will never
create a forbidden subconfiguration. Suppose we had zero probability
of seeing any $\Cmax$ rungs. Consider the measure
$\mu'$ obtained by changing each rung to $\Cmax$ independently
with some small probability $0 < \eps < 1$. Then $\mu'$ is also 
ergodic, and a straightforward computation shows that its
measure theoretic entropy is 
$h(\mu') = (1 - \eps) h(\mu) + H(\eps)$, where 
$H(\eps) = - \eps \log \eps - (1 - \eps) \log (1 - \eps)$.
Hence for $\eps$ sufficiently small, $h(\mu') > h(\mu)$,
a contradiction.

Consider now the sequence of ``blocks" between successive
$\Cmax$ rungs. These form a stationary sequence. Again, by 
maximum entropy, the blocks have to be independent. Indeed, if
they were not, consider the measure $\mu'$, where the blocks 
are i.i.d.\ and each block has its $\mu$-distribution. Let 
$E := \{ C_0 = \Cmax \}$. Since the expected
length of a block is the same in $\mu'$ and $\mu$, we have
$\mu'(E) = \mu(E)$. The
measure-preserving maps induced by $E$ \cite[Chapter 1, \S 5]{CFS82}
are the translations of blocks, with invariant measures
$\mu_E$ and $\mu'_E$ (the normalized restrictions of $\mu$ and $\mu'$
to $E$). Since $\mu'_E$ is i.i.d.~and
$\mu_E$ is not, we get
\eqnst
{ h(\mu') 
  = \mu'(E) h(\mu'_E) 
  > \mu(E) h(\mu_E) 
  = h(\mu) }
by \cite[Chapter 10, \S 6, Theorem 2]{CFS82}.

It follows that $\mu$ is determined by the joint distribution 
of renewal times (distance between $\Cmax$ rungs) and the
inter-renewal configuration. Suppose that a block has positive
probability of being non-left-burnable and also positive 
probablity of being non-right-burnable. Then with probability 
one, there will be a non-left-burnable block to the left of
a non-right-burnable block. This creates a forbidden subconfiguration,
and hence is
impossible. Therefore, at most one of the above possibilities 
has positive probability. Assume without loss of generality
that blocks are left-burnable with probability $1$.

Consider now the configuration between two renewals (not 
necessarily consecutive) that are distance $L$ apart. By maximum
entropy, the conditional distribution of the configuration
given $L$ is uniform over all left-burnable 
configurations of length $L-1$. Since this holds for
arbitrarily large $L$, it implies that the finite-dimensional 
distributions of $\mu$ are given by the thermodynamic limit
of $\muL_{n,m}$, and hence $\mu = \muL$. Analogously, we get
$\mu = \muR$ if blocks are right-burnable with probability $1$.

The proof in the symmetric case is very similar. Adding $\Cmax$
rungs in an i.i.d.~fashion does not destroy the symmetry of the
measure, and hence $\Cmax$ rungs have to occur with positive 
probability. As before, they are renewals. Again, blocks have to
be either left- or right-burnable, and by symmetry, they have to
be both with probability $1$. As before, this implies that 
the measure coincides with $\muS$.
\end{proof}

\begin{theorem}
\label{thm:therm}
If $-n,m \to \infty$ in such a way that $\lim -n/m = \rho / (1-\rho)$,
$\rho \in [0,1]$, then 
\eqnst
{ \lim \mu_{n,m}
  = \rho \muL + (1 - \rho) \muR. }
Consequently, the set of weak limit points of $\{ \mu_{n,m} \}$ 
consists of all convex combinations of 
$\muL$ and $\muR$. 
\end{theorem}

\begin{proof}
The idea of the proof is the following. A recurrent configuration
has to burn if we burn from both the left and the right. We show that
the left- and right-burnable portions of the configuration almost form
a partition of $\Lam_{n,m}$, up to an overlap or uncovered region of 
size $o(m-n)$ in probability, and that the location of the ``boundary 
layer" between them is approximately uniform over $[n,m]$. This implies 
that in a fixed finite window we see a convex combination of $\muL$ 
and $\muR$.

For $\eta \in \Omega_{n,m}$, let
\eqnsplst
{ \sigmaL 
  &:= \sigmaL_{n,m}
  := \max \{ k : \text{$C_k = \Cmax$ and 
    $\eta_{\Lam_{n,k}}$ is left-burnable} \}, \\
  \sigmaR
  &:= \sigmaR_{n,m}
  := \min \{ k : \text{$C_k = \Cmax$ and 
    $\eta_{\Lam_{k,m}}$ is right-burnable} \}, }
where the values $n - 1$ and $m + 1$ are allowed in both 
cases if there is no $k$ with the required property.

We show that 
\eqn{e:narrow}
{ \frac{|\sigmaL - \sigmaR|}{m - n} \to 0 }
in probability.

\emph{Case 1: $\sigmaL < \sigmaR$.}
We show that the number of possible configurations between 
$s := \sigmaL$ and $t := \sigmaR$ has exponential growth rate
smaller than $h^{\rm L}$. 

Let $\eta^0 := \eta_{\Lam_{s+1,t-1}}$. Consider 
left-burning on $\eta^0$. By the definition of $\sigmaL$, and
since $C_t = \Cmax$, the rightmost site of 
$\eta^0$ that will be left-burnt (when left-burning $\eta^0$)
is in a rung $k$ with
$s \le k < t - 1$. Here $k = s$ if left-burning cannot
start. Similarly, the leftmost site
of $\eta^0$ that can be right-burnt is in a rung $l$ with
$s + 1 < l \le t$. Since $\eta^0$ is burnable,
we need to have $l \le k+1$. 

Let $\xi^0$ be the configuration obtained by replacing rung 
$k + 1$ of $\eta^0$ by $\Cmax$. Since $\xi^0$ is also burnable, 
it follows easily that $\eta^0_{\Lam_{s+1,k}}$ is
left-burnable, and $\eta^0_{\Lam_{k+2,t-1}}$ is 
right-burnable. Therefore, the number of possibilities for
$\eta^0$ is bounded by
\eqnst
{ \sum_{k = s}^{t-2} b_{k-s} |\cC| b_{t-k-2}
  \le C e^{(h^{\rm L} - \delta) (t-s) } }
for some $\delta > 0$ and some $C < \infty$ by Lemma \ref{lem:entropy-compare}.

Summing over all possible values of $s$ and $t$, it follows that 
for any $\eps > 0$ there exist $C_1 = C_1 (\eps)$ 
and $c_1 = c_1 (\eps) > 0$ such that
\eqn{e:missing}
{ \mu_{n,m} \{ \sigmaR - \sigmaL \ge \eps (m - n) \} 
  \le C_1 e^{- c_1 (m - n)}. }

\emph{Case 2: $\sigmaL \ge \sigmaR$.} Observe that the
configuration between $\sigmaR$ and $\sigmaL$ is both
left-burnable and right-burnable, hence it belongs to
$\OmegaS_{\sigmaR,\sigmaL}$. Also, the configuration
to the left of $\sigmaR$ is in $\OmegaL_{n,\sigmaR}$,
and the configuration to the right of $\sigmaL$ is in
$\OmegaR_{\sigmaL,m}$. 

Since $h^{\rm S} < h^{\rm L}$, it follows that for any $\eps > 0$ there 
exists $C_2$ and $c_2 = c_2(\eps) > 0$ such that
\eqn{e:overlap}
{ \mu_{n,m} \{ \sigmaL - \sigmaR \ge \eps (m - n) \} 
  \le C_2 e^{-c_2(m - n)}. }

The bounds \eqref{e:missing} and \eqref{e:overlap} establish
\eqref{e:narrow}.

In the remainder of the proof we are going to need a minor
variation on \eqref{e:overlap} when $\sigmaR \le \sigmaL$.
The reason is that the value of $\sigmaR$ gives some 
information on the left-burnable configuration to the left of 
$\sigmaR$ (namely, that it is not right-burnable if it contains a rung
$\Cmax$), whereas we
would like to achieve independence. 
Let $\hat{\sigma}^{\rm R}$ denote the rightmost $\Cmax$ rung to the left 
of $\sigmaR$ (we set $\hat{\sigma}^{\rm R} = n - 1$ if such a rung 
does not exist). Then the configuration between $\hat{\sigma}^{\rm R}$ 
and $\sigmaR$ is left-burnable but not right-burnable. 
In any case, it is in
$\Omega^{\rm L,0}_{\hat{\sigma}^{\rm R}+1,\sigmaR-1}$. We define
$\hat{\sigma}^{\rm L}$ analogously. 
By similar arguments as before, we have the bound
\eqn{e:overlap2}
{ \mu_{n,m} \{\sigmaL \ge \sigmaR \hbox{ and }
  \hat{\sigma}^{\rm L} - \hat{\sigma}^{\rm R} \ge \eps (m - n) \} 
  \le C_3 e^{ -c_3 (m - n)}. }

Next we prove that the location of the ``boundary layer" 
between the left- and right-burnable parts is approximately 
uniform.

First condition on the value of $d := \sigmaR - \sigmaL$ in the 
case when $d$ is positive. Observe that given $\sigmaL = s$ and
$\sigmaR = t$, the configurations on $\Lam_{n,s}$,
$\Lam_{s,t}$ and $\Lam_{t,m}$ are conditionally independent. 
Also, the configuration on
$\Lam_{n,s-1}$ has law $\muL_{n,s-1}$ and the configuration
on $\Lam_{t+1,m}$ has law $\muR_{t+1,m}$. Noting that 
$\muR$ is the reflection of $\muL$, we can
uniquely represent the configuration 
in the following way. Draw a sample $\eta$ from
$\muL_{n,m-d}$ conditioned on having at least one renewal.
Select one of the $\Cmax$ rungs uniformly at random: suppose it is 
rung $S$. Draw an independent sample $\xi$ from the set of 
configurations $\eta^0$ described under Case 1 above having length $d-1$. 
Concatenate the configurations $\eta_{\Lam_{n,S}}$, $\xi$, $\Cmax$, 
and the reversal of $\eta_{\Lam_{S+1,m-d}}$. This gives all
configurations with $\sigmaL = S$ and $\sigmaR = S + d$, and
the representation is unique.

Next we want to show that the random variable $S$ defined above is 
roughly uniformly distributed in $[n,m-d]$. First note that by 
Lemma \ref{lem:renewals-finite}, under $\muL_{n,m-d}$,
the distribution of the sequence of inter-renewal times is exchangeable.
Also, due to the inequality $h^{\rm L,0} < h^{\rm L}$, the longest
inter-renewal time is $o(m - n - d)$ in probability.
These two together imply that for any $0 < u < 1$, 
\eqnst
{ \muL_{n,m-d} \{ S-n < u (m - n - d) \} 
  \to u \quad \text{as $m - n \to \infty$} }
uniformly in $1 \le d < (m-n)/2$.
This implies that 
\eqnst
{ \mu_{n,m} \{ \sigmaL-n < u (m - n - d) \,|\, \sigmaR - \sigmaL = d \}
  \to u \quad \text{as $m - n \to \infty$} }
uniformly in $1 \le d < (m-n)/2$. Averaging over $1 \le d \le \eps(m-n)$,
we get
\eqn{e:uniform1}
{ \mu_{n,m} \left\{ \frac{\sigmaL-n}{m - n} < u \,\Biggm|\, 
    1 \le \sigmaR - \sigmaL \le \eps(m-n) \right\}
  = u + O(\eps) + o(1) }
as $m - n \to \infty$. 

Now condition on $d := \hat{\sigma}^{\rm L} - \hat{\sigma}^{\rm R}$ in the
case when $\sigmaL \ge \sigmaR$. Given 
$\hat{\sigma}^{\rm R} = \hat{s}^{\rm R}$ and $\hat{\sigma}^{\rm L} = \hat{s}^{\rm L}$,
the configurations on $\Lam_{n,\hat{s}^{\rm R}-1}$, 
$\Lam_{\hat{s}^{\rm R},\hat{s}^{\rm L}}$ and $\Lam_{\hat{s}^{\rm L}+1,m}$ are conditionally
independent, with the first and the third having laws
$\muL_{n,\hat{s}^{\rm R}-1}$ and $\muR_{\hat{s}^{\rm L}+1,m}$ (respectively).
Therefore, the configuration can be represented analogously to the
case $\sigmaL < \sigmaR$, which gives rise to the estimate
\eqn{e:uniform2}
{ \mu_{n,m} \left\{ \frac{\hat{\sigma}^{\rm R}-n}{m - n} < u \,\Biggm|\, 
    \hat{\sigma}^{\rm L} - \hat{\sigma}^{\rm R} \le \eps(m-n) \right\}
  = u + O(\eps) + o(1) }
as $m - n \to \infty$.

We are ready to complete the proof of the theorem. Suppose we
have a cylinder event $E$ depending on the configuration in 
$\Lam_{-k,k}$. Let
\eqnst
{ \tau^{\rm L}
  := \begin{cases}
    \sigmaL & \text{if $\sigmaL < \sigmaR$}, \\
    \hat{\sigma}^{\rm R} & \text{if $\sigmaL \ge \sigmaR$},
  \end{cases}
  \qquad\mbox{and}\qquad
 \tau^{\rm R}
  := \begin{cases}
    \sigmaR & \text{if $\sigmaL < \sigmaR$}, \\
    \hat{\sigma}^{\rm L} & \text{if $\sigmaL \ge \sigmaR$}.
  \end{cases} }
Let $A_\eps := \{ \tau^{\rm L} > \eps(m-n) \}$ and
$B_\eps := \{ \tau^{\rm R} < -\eps(m-n) \}$. For $t > \eps(m-n) > k$,
\eqnst
{ \mu_{n,m} \{ E \,|\, A_\eps,\, \tau^{\rm L} = t \}
  = \muL_{n,t-1} \{ E \}
  = \muL \{ E \} \big(1 + o_\eps(1)\big) }
as $n \to -\infty$ and $m \to \infty$, where the $o_\eps(1)$ depends 
on $\eps$, but not on $t$. Similarly, for $t < -\eps(m-n) < -k$,
\eqnst
{ \mu_{n,m} \{ E \,|\, B_\eps,\, \tau^{\rm R} = t \}
  = \muR_{t+1,m} \{ E \}
  = \muR \{ E \} \big(1 + o_\eps(1)\big). }
Since by \eqref{e:narrow}, \eqref{e:overlap2},
\eqref{e:uniform1} and \eqref{e:uniform2}, 
$\mu_{n,m} \{ A_\eps \} = -n / (m-n) + O(\eps)$ and
$\mu_{n,m} \{ B_\eps \} = m / (m-n) + O(\eps)$, the
theorem follows by letting
$\eps \to 0$.
\end{proof}

\section{Avalanches}
\label{sec:avalanches}

By toppling in an infinite graph we mean the following. 
Suppose we start from a configuration $\eta$ with finitely
many unstable sites. We simultaneously topple all unstable 
sites, and repeat this as long as there are unstable sites
(possibly infinitely many times). After each step there are
only finitely many unstable sites. This is equivalent to
toppling sites one-by-one during each step, before moving
on to toppling other sites. Let us call this the 
\emph{standard toppling}.

\begin{defn}
A (possibly infinite) sequence of topplings is called
\emph{legal} if it has the properties: (i) only unstable
sites are toppled in each step; (ii) any site that is unstable 
at some step will be toppled at some later step. 
\end{defn}

\begin{lemma}
Any two legal sequences of topplings are equivalent in the
sense that each site topples the same number of times in
both sequences (which may be infinity). In particular,
any legal sequence of topplings is equivalent to 
standard toppling.
\end{lemma}

\begin{proof}
This can be proved the same way as for finite sequences 
of topplings \cite{MRZ02}. Given two legal sequences of 
topplings at the sites
\eqnsplst
{ &x_1, x_2, \dots \\
  &y_1, y_2, \dots }
we can transform one into the other. Since $x_1$ is unstable
at the beginning, it has to occur in the second sequence.
Suppose it occurs first as $y_{k_1}$. Then the toppling of
$y_{k_1}$ can be commuted through the topplings of 
$y_1, y_2, \dots, y_{k_1 - 1}$, so the $y$-sequence is
equivalent to 
\eqnst
{ x_1, y_2, \dots, y_{k_1-1}, y_{k_1+1}, y_{k_1+2}, \dots }
We can now eliminate $x_1$ from both sequences, and the
lemma follows.
\end{proof}

\begin{theorem}
Suppose we add one grain to each site in rung $0$ in an infinite
left-burnable configuration. Then each site will topple infinitely 
many times. The same holds for right-burnable configurations.
\end{theorem}

\begin{proof}
Add one grain to each site in rung $0$, and initially, 
do not topple in rungs to the left of zero. (For the 
moment, let us disregard that this may be an illegal sequence
of topplings.) The topplings that occur 
on the right are equivalent to the burning procedure on
$\Lam_{0,\infty}$. Since the configuration is left-burnable,
each site in $\Lam_{0,\infty}$ will topple exactly once.
In particular, each site in rung $0$ will have toppled.
Also, it is easy to verify that each site in $\Lam_{0,\infty}$ 
will have received as many grains as it has lost, and
hence has its original height.

The topplings in rung $0$ give one grain to each site
in rung $-1$. Therefore, the argument can be repeated
as if we have added one grain to each site in rung $-1$,
and hence topplings continue forever.
This almost completes the argument, apart from the
technicality that this is not a legal sequence of topplings.
Instead, now we carry out the topplings on the right only
to a large finite time until a rung $p_1 \gg 1$ is toppled.
Then carry out topplings started from rung $-1$, until
a rung $1 < p_2 < p_1$ is toppled, and so on. If $p_1 \ge 2K$,
we can repeat this with
rungs $1 < K < p_K < \dots < p_1$, for any given large $K$. 
At this point, rung $-k$ has toppled $K-k$ times for
$k = 1, \dots, K$, and rungs $0, 1, \dots, K$ have toppled 
$K$ times. It follows that in any legal sequence of topplings,
sites $-K/2,\dots,K/2$ each topple at least $K/2$ times. Since
$K$ was arbitrary, the theorem follows.
\end{proof}

\begin{remark}
As the following example shows, there can be infinite avalanches
such that every site topples only finitely many times. Take 
$G = I_{0,1}$. Under $\mu^L$, there is positive probability that the 
configuration at rungs $1$--$6$ equals
\eqnst
{ \begin{matrix}
  3 & 2 & 3 & 3 & 1 & 3 \\
  3 & 3 & 1 & 3 & 3 & 3  
  \end{matrix} }
Now adding a grain to the first row in rung $4$ yields an avalanche 
with toppling numbers:
\eqnst
{ \begin{matrix}
  \dots & 0 & 0 & 1 & 2 & 1 & 1 & \dots \\
  \dots & 0 & 0 & 0 & 1 & 1 & 1 & \dots  
  \end{matrix} }
\end{remark}

\section{Coding by Markov chains}
\label{sec:coding}

In this section we show that the measures $\muL_{n,m}$ and
$\muL$ can be coded by a Markov chain with finitely many
states. Before proving this for a general graph $G$,
we sketch a proof in the special case $G = I_{0,1}$. Although 
for general $G$ we will not have as explicit a description as 
for $I_{0, 1}$, the approach will be similar.

\medbreak

\emph{Coding by a finite Markov chain for $G = I_{0,1}$.}
Based on Lemma \ref{lem:left-burn}, the following equivalent
description of left-burnable configurations can be given. 
Consider the alphabet of symbols 
\eqn{e:states}
{ \cA 
  := \{ (3,3), (3,2), (2,3) , (3,1), (1,3), \overline{(3,2)}, 
    \overline{(2,3)} \}\,. } 
Let $\cA_{n,m} := \cA^{I_{n,m}}$.
We think of $\overline{(3,2)}$ replacing a $(3,2)$ rung
that is following a $(3,1)$ before the next $(3,3)$ occurs.
It follows from the characterization in Lemma \ref{lem:left-burn}
that elements of $\OmegaL_{n,m}$ can be coded in a one-to-one
fashion by a set $\bOmega^{\rm L}_{n,m} \subset \cA_{n,m}$ that is 
a topological Markov chain (subshift of finite type) 
\cite[Section 1.9]{KH95} with alphabet $\cA$. Namely, the only 
restrictions on sequences in $\bOmega^{\rm L}_{n,m}$ are that certain 
pairs of symbols cannot occur next to each other. For example:
(a) $(3,3)$ has to be followed by $(3,3)$, $(3,2)$, $(2,3)$, $(3,1)$
or $(1,3)$; (b) $(3,1)$ has to be followed by $(3,3)$ or
$\overline{(3,2)}$; (c) $\overline{(3,2)}$ has to be followed by
$\overline{(3,2)}$ or $(3,3)$; etc. The full transition matrix is
\eqn{e:matrix}
{ T := \left( \begin{matrix}
  1 & 1 & 1 & 1 & 1 & 0 & 0 \\
  1 & 1 & 1 & 1 & 1 & 0 & 0 \\
  1 & 1 & 1 & 1 & 1 & 0 & 0 \\
  1 & 0 & 0 & 0 & 0 & 1 & 0 \\
  1 & 0 & 0 & 0 & 0 & 0 & 1 \\
  1 & 0 & 0 & 0 & 0 & 1 & 0 \\
  1 & 0 & 0 & 0 & 0 & 0 & 1 
\end{matrix} \right), }
where the rows and columns correspond to the symbols in the order
displayed in \eqref{e:states}. Due to the special role of the 
symbols $\overline{(3,2)}$ and $\overline{(2,3)}$, we need to add 
the boundary condition that the rung at $n$ is not one of these.

It is not hard to check that the topological Markov chain is
transitive \cite[Definition 1.9.6]{KH95}; in fact, all entries of 
$T^3$ are positive. 
Let $\bar{\Omega}^{\rm L} \subset \cA^\Z$ denote the subshift
defined by $T$, and let $\bar{\mu}^{\rm L}$ be its Parry measure, which is
a Markov chain. 
By \cite[Section 4.4]{KH95}, $\bar{\mu}^{\rm L}$ is the unique measure 
of maximum entropy on $\bar{\Omega}^{\rm L}$. Let 
${\cal P} : \bar{\Omega}^{\rm L} \to \OmegaL$ denote the map that
replaces each $\overline{(3,2)}$ by $(3,2)$ and each 
$\overline{(2,3)}$ by $(2,3)$. Since $\muL$ has maximal entropy 
by Theorem \ref{thm:max-entropy}, and ${\cal P}^{-1}$ is well
defined $\muL$-almost surely, ${\cal P}$ is a metric isomorphism 
between $\bar{\mu}^{\rm L}$ and $\muL$.

\medbreak

Now we generalize the coding to an arbitrary graph $G$. 
First note that it is not very surprising that such a coding
should exist. Using Majumdar and Dhar's tree construction \cite{MD92},
recurrent configurations in $\Lam_{n,m}$ are in one-to-one
correspondence with spanning trees of $\Lam_{n,m}$ with
wired boundary conditions. It has been shown in \cite{OHthesis}
that spanning trees have a Markovian coding.
However, since the correspondence is non-local, it does not seem 
easy to deduce a Markovian coding from the spanning-tree result.

\medbreak

We let $\cP := \cP(G)$ denote the set of all subsets of $G$.

\begin{theorem}
\label{thm:coding}
There exists an alphabet $\cA := \cA(G) \subset \cC \times \cP \times
\cP^\cP$, 
an inclusion $i: \cC \to \cC \times \cP \times \cP^\cP$, and a 
transitive $0$-$1$ matrix $T := T(G)$ indexed by $\cA$ such that 
for each $m$, the set
$\OmegaL_{1,m}$ is in one-to-one correspondence with the set of 
sequences 
\eqnsplst
{ \bOmega^{\rm L}_{1,m}
  &:= \bOmega^{\rm L}_{1,m}(G) \\
  &:= \{ \omega \in \cA_{1,m} :
    \omega_1 \in i(\cC),\,
    T(\omega_k, \omega_{k+1}) = 1,\,
    k = 1, \dots, m-1 \}. }
The correspondence is given by the projection
$P : \cC \times \cP \times \cP^\cP \to \cC$ applied coordinatewise.
\end{theorem}

For the proof of Theorem \ref{thm:coding}, we will need to perform 
left-burning in a special way, as introduced below. This can be 
regarded as a generalization of the rung-by-rung argument from 
the proof of Lemma \ref{lem:left-burn}. Following the definition of
the special burning rule, we use it to prove two lemmas that will
lead to the proof of Theorem \ref{thm:coding}.
Once Theorem \ref{thm:coding} is established, the Markov chain that codes
$\muL$ is the Parry measure, as for $G = I_{0, 1}$.

\medbreak

\emph{Burning with leftmost rung rule.}
We perform burning one rung at a time, with the
rule that whenever there are no more burnable sites in the rung
currently being burnt, we move on to the leftmost rung that
has burnable sites. We now describe the procedure in more detail.

We first burn sites in rung $1$ that can be burnt consistent
with the left-burning rule. When there are no more burnable sites
in rung 1, we start burning sites in rung 2, and continue 
burning rung 2 until there are no more burnable sites in that
rung. This may have created further burnable sites in
rung 1. If there are such, we burn sites in rung 1, again until 
there are no more burnable sites in that rung. At some point
there will be no burnable sites in either rung 1 or 2. Now
we burn sites in rung 3, and move between rungs 1, 2 and 3
until there are no more burnable sites in those rungs.
In general, we move on to rung $k+1$ when there are no more
burnable sites in rungs $1, \dots, k$.

If the configuration we started with is $C_1 \vee \dots \vee C_m$,
we adopt the following convention for burning the rightmost
rung $C_m$. We add a ``ghost" rung $C_{m+1} = \Cmax$
that will remain unburnt until the first time there are no
more burnable sites in rungs $1, \dots, m$. At this time, 
we burn the ghost rung, and continue with the leftmost
rule. It is easy to see that this yields an equivalent definition
of left-burnability, that is, all rungs will burn if and only if
the original configuration was left-burnable.

For $1 \le k \le m$ and $C_1 \vee \dots \vee C_k$ left-burnable,
let $T_k + 1$ be the first time we burn a site in rung $k+1$. 
It is easy to see that all rungs are burnt at
time $T_{m+1}$ if and only if $\eta \in \OmegaL_{1,m}$.

\medbreak

Before stating the two lemmas needed for Theorem \ref{thm:coding},
we need some notation.
Let $\eta = C_1 \vee \dots \vee C_m$ be a configuration with
$C_k \in \cC$, $1 \le k \le m$. Let $B_k \subset G$ denote the set of 
sites in rung $k$ that have been burnt by time $T_k$. 
The sequence $(C_k,B_k)_{k=1}^m$ is non-Markovian in general.
We note, however, in order to motivate the arguments to come, that 
if $G = I_{0,1}$, then $B_{k+1}$ is a function of
$C_k$, $C_{k+1}$ and $B_k$ only (it depends on 
$C_1, \dots, C_{k-1}$ only through $C_k$, $C_{k+1}$ and $B_k$). 
It is not hard to show that this implies that $(C_k,B_k)_{k=1}^m$ 
is Markovian. The proof is similar (and simpler) than that of 
Theorem \ref{thm:coding} below, and is left to the reader. 
For general $G$, our strategy will be to augment the information 
contained in $B_k$ so that we get a Markovian sequence.

Fix $(C_j, B_j)_{j = 1}^k$, where $1 \le k \le m$. Depending on 
this sequence, we define a function $f_k : \cP \to \cP$
that will encode what the effect is of burning in rung $k + 1$
on the future of the burning process in rungs $1 \le j \le k$. 
We stress that the definition of $f_k$ will
ignore the actual value of $C_{k+1}$; in particular, it will
also make sense for $k = m$. Fix $A \subset G$. 
Regardless of the value of $C_{k+1}$, let us declare all sites 
in $A \times \{ k + 1 \}$ to be burnt. This may create burnable
sites in rung $k$ after $\bigcup_{j=1}^k B_j$ has been burnt.
Now let us perform burning with the
leftmost rung rule until there are no more burnable sites
in rungs $1 \le j \le k$. This process does not use information 
about rung $k + 1$ other than the specified set $A$.
We define $f_k(A)$ to be the set of sites that are burnt 
in rung $k$ at the end of this process. For example, we  
have $f_k(\es) = B_k$, and more generally, $f_k(A) = B_k$ for
$A \subset B_k$, since in this case no new burnable 
sites appear in rung $k$. Whenever $C_1 \vee \dots \vee C_k$ 
is left-burnable, we have $f_k(G) = G$. In general, we have 
$B_k \subset f_k(A) \subset G$ for $A \in \cP$.

We prove Theorem \ref{thm:coding} by showing that 
$(C_k,B_k,f_k)_{k=1}^m$ is Markovian. We verify this in the two
lemmas below that characterize the pairs that can occur next to 
each other for left-burnable $\eta = C_1 \vee \dots \vee C_m$. 
To facilitate the proof, we define an auxiliary function 
$g : \cP \times \cC \times \cP \to \cP$. Given $A, A' \subset G$ 
and $C \in \cC$, we set the configuration in rung $1$ to be $C$ 
and declare all sites in $A \times \{ 0 \} \cup A' \times \{ 2 \}$ 
to be burnt. Now we perform left-burning in rung $1$. By this we 
mean specifying a maximal sequence of vertices 
$v_1, \dots, v_k \in G \times \{ 1 \}$, such that 
the requirements of Definition \ref{defn:left-burn} are satisfied with 
$V := (G \setminus A) \times \{ 0 \} \cup G \times \{ 1 \} \cup
  (G \setminus A') \times \{ 2 \}$. We define
$g(A, C, A')$ to be the set of sites that burn in rung $1$.

\begin{lemma}
\label{lem:whatpairs}
For $\eta = C_1 \vee \dots \vee C_m
\in \OmegaL_{1,m}$, the following properties hold:
\begin{itemize}
\item[(a)] $(B_1, f_1) = \psi(C_1)$ for some function $\psi = \psi_G$,
  in fact, $B_1 = g(G,C_1,\es)$ and $f_1(A) = g(G,C_1,A)$;
\item[(b)] $g(B_k,C_{k+1},\es) \not= \es$, $1 \le k < m$; 
\item[(c)] $(B_{k+1}, f_{k+1}) = \phi (B_k, C_{k+1}, f_k)$ for a 
  function $\phi = \phi_G$ independent of $k$, $1 \le k < m$; and
\item[(d)] $f_k(G) = G$, $1 \le k \le m$.
\end{itemize}
\end{lemma}

\begin{proof}
(a) follows directly from the definitions, and (d) has been
observed before the statement of the lemma.
If (b) failed for some $1 \le k < m$, that 
would mean that after time $T_k$ there were no burnable sites
in rungs $1, \dots, k+1$, with rung $k+1$ completely
unburnt. That means that there are no burnable sites at all
after time $T_k$, which contradicts the burnability of $\eta$.

The proof of (c) is a bit lengthy. We first show that $B_{k+1}$ is
a function of $B_k$, $C_{k+1}$ and $f_k$. 
For this, we look at the burning process between times $T_k + 1$ and
$T_{k+1}$ in more detail. We define the following intermediate 
times: letting $R_0 := T_k + 1$, we define
\eqnsplst
{ S_1 
  &:= \min \left\{ R_0 \le n \le T_{k+1} :\ \parbox{2in}
     {there are no burnable sites in rung $k+1$ at time $n$} 
     \right\}; \\
  R_1
  &:= \min \left\{ S_1 \le n \le T_{k+1} :\ \parbox{2.1in}
     {there are no burnable sites in rungs $1 \le j \le k$ 
     at time $n$} \right\}, }
and for $i \ge 2$ we recursively set
\eqnsplst
{ S_i 
  &:= \min \left\{ R_{i-1} \le n \le T_{k+1} :\ \parbox{2in}
     {there are no burnable sites in rung $k+1$ at time $n$} 
     \right\}; \\
  R_i
  &:= \min \left\{ S_i \le n \le T_{k+1} :\ \parbox{2.1in}
     {there are no burnable sites in rungs $1 \le j \le k$ 
     at time $n$} \right\}. }
Set $B^{(0)} := B_k$. Between times $T_k + 1 = R_0$ and $S_1$, 
the subset $A^{(0)} := g(B^{(0)}, C_{k+1}, \es)$ of rung 
$k+1$ is burnt. Let $B^{(1)} := f_k(A^{(0)}) \supset B^{(0)}$. 
By the definition of $f_k$, $B^{(1)}$ is the set of sites in
rung $k$ that is burnt by time $R_1$. We set 
$A^{(1)} := g(B^{(1)}, C_{k+1}, \es)$ and $B^{(2)} := f_k(A^{(1)})$. 
By the definition of $g$, $A^{(1)}$ is the set of sites in rung 
$k+1$ burnt at time $S_2$. Although less obvious, $B^{(2)}$ is
the set of sites in rung $k$ burnt at time $R_2$. The latter
statement needs careful proof since $f_k$ was defined in terms 
of the state of the burning process at time $R_0$ rather than at
$R_1$. Consider the sequence of sites burnt in the computation 
of $f_k(A^{(1)})$ (following the definition). We merely get a 
rearrangement of this sequence if we first declare 
$A^{(0)} \times \{ k + 1 \}$ to be burnt, let burning act on 
rungs $1, \dots, k$, then declare 
$\left( A^{(1)} \setminus A^{(0)} \right) \times \{ k + 1 \}$
to be burnt, and then let burning act on rungs $1, \dots, k$.
This observation proves our claim about $B^{(2)}$.

In general, for $i \ge 1$, after burning between times 
$R_i$ and $S_{i+1}$, the set of sites burnt in rung $k+1$ is 
$A^{(i)} := g( B^{(i)}, C_{k+1}, \es) \supset A^{(i-1)}$. 
We set $B^{(i+1)} := f_k(A^{(i)})$. Similarly to the case
$i = 1$ spelled out above, by a decomposition of $A^{(i)}$,
we get that $B^{(i+1)}$ is the set of sites burnt in rung $k$ 
at time $R_{i+1}$. Since there is some $j_0$ for which 
$B^{(j+1)} = B^{(j)}$ and
$A^{(j+1)} = A^{(j)}$ for $j \ge j_0$, we have 
$B_{k+1} = A^{(j_0)}$. To summarize, $B_{k+1}$ is obtained 
as the stable result of applying the functions $g$ and $f_k$ 
according to 
\eqn{e:Bk+1}
{\begin{aligned}
  B^{(0)} &:= B_k\,, & 
  A^{(0)} &:= g(B^{(0)}, C_{k+1}, \es)\,, \\
  B^{(1)} &:= f_k(A^{(0)})\,, &
  A^{(1)} &:= g(B^{(1)},C_{k+1}, \es)\,, \\
  &\vdots & &\vdots \\
  B^{(j_0)} &:= f_k(A^{(j_0-1)})\,, & 
  B_{k+1} &= A^{(j_0)} := g(B^{(j_0)}, C_{k+1}, \es)\,. 
  \end{aligned} }
This shows that $B_{k+1}$ is a function of 
$B_k$, $C_{k+1}$ and $f_k$ that does not depend on the value of $k$ (is
``$k$-independent") .

Now we can prove the remainder of
(c) by a similar argument. Consider the state of
the burning process at time $T_{k+1}$, at which time the
set of sites burnt in rung $k + 1$ is $B_{k+1}$.
By the defintion of $f_k$, the set of sites in rung $k$
burnt at time $T_{k+1}$ is the set $\bar{B}^{(0)} := f_k(B_{k+1})$.
Given $A \subset G$, declare all sites in $A \times \{k + 2\}$ 
to be burnt (ignoring $C_{k+2}$). Now perform burning
in rung $k+1$, which ends at some time $\bar{S}_1$. 
Then the set of sites in rung $k+1$ burnt at time 
$\bar{S}_1$ is 
\eqnst
{ \bar{A}^{(0)} 
  := g(\bar{B}^{(0)}, C_{k+1}, A) 
    \supset B_{k+1}. }
Now we can essentially apply the argument above starting with 
$\bar{A}^{(0)}$ in place of $A^{(0)}$.
We perform burning on rungs $1, \dots, k$, and let 
$\bar{R}_1 \ge \bar{S}_1$ be the first time when there
are no burnable sites in these rungs.
Then the set of sites in rung $k$ burnt at time $\bar{R}_1$ is 
$\bar{B}^{(1)} := f_k(\bar{A}^{(0)})$. This is shown by
the decomposition
\eqnst
{ \bar{A}^{(0)}
  = A^{(0)} \cup 
    \left[ \bigcup_{j = 1}^\infty 
    \left( A^{(j)} \setminus A^{(j-1)} \right) \right]
    \cup \left( \bar{A}^{(0)} \setminus B_{k+1} \right). }
Next we perform burning in rung $k+1$ that stops at some
time $\bar{S}_2 \ge \bar{R}_1$, and then on rungs 
$1, \dots, k$, which stops at $\bar{R}_2 \ge \bar{S}_2$. 
The set of sites in rung $k+1$ burnt at time $\bar{S}_2$ is 
$\bar{A}^{(1)} := g(\bar{B}^{(1)}, C_{k+1}, A)$, and the
set of sites in rung $k$ burnt at time $\bar{R}_2$ is
$f_k(\bar{A}^{(1)})$. We continue to iterate $g$ and $f_k$
until the burnt sites in rung $k+1$ stabilize to some 
set $\bar{A}^{(\bar{j}_0)}$. We then have 
$f_{k+1}(A) = \bar{A}^{(\bar{j}_0)}$. We have
\eqn{e:fk+1}
{\begin{aligned}
  \bar{B}^{(0)} &:= f_k(B_{k+1})\,, & 
  \bar{A}^{(0)} &:= g(\bar{B}^{(0)}, C_{k+1}, A)\,, \\
  \bar{B}^{(1)} &:= f_k(\bar{A}^{(0)})\,, &
  \bar{A}^{(1)} &:= g(\bar{B}^{(1)},C_{k+1}, A)\,, \\
  &\vdots & &\vdots \\
  \bar{B}^{(\bar{j}_0)} &:= f_k(\bar{A}^{(\bar{j}_0-1)})\,, & 
  f_{k+1}(A) &= \bar{A}^{(\bar{j}_0)} 
    := g(\bar{B}^{(\bar{j}_0)}, C_{k+1}, A)\,. 
  \end{aligned} }
This shows that $f_{k+1}$ is a $k$-independent function 
of $B_k$, $C_{k+1}$ and $f_k$, and hence (c) follows.
\end{proof}

\begin{lemma}
\label{lem:burnable}
Let $\psi$ and $\phi$ be as in Lemma \ref{lem:whatpairs}.
Suppose the sequence $(C'_k,B'_k,f'_k) \in \cC \times \cP \times \cP^\cP$,
$k = 1, \dots, m$, satisfies the conditions:
\begin{itemize}
\item[(A)] $(B'_1,f'_1) = \psi(C'_1)$;
\item[(B)] $g(B'_k,C'_{k+1},\es) \not= \es$, 
  $1 \le k < m$; 
\item[(C)] $(B'_{k+1},f'_{k+1}) = \phi(B'_k,C'_{k+1},f'_k)$,
  $1 \le k < m$; and
\item[(D)] $f'_k(G) = G$, $1 \le k \le m$.
\end{itemize}
Then $\eta := C'_1 \vee \dots \vee C'_m \in \OmegaL_{1,m}$,
and taking $C_k := C'_k$ in the definitions preceding
Lemma \ref{lem:whatpairs}, we have $B_k = B'_k$ and $f_k = f'_k$,
$1 \le k \le m$.
\end{lemma}

\begin{proof}
We verify the statement by induction on $m$. When $m = 1$,
$\eta = C'_1 \in \OmegaL_{1,1}$ since $C'_1 \in \cC$.
Therefore by (A), $(B_1,f_1) = \psi(C_1) = \psi(C'_1)
= (B'_1,f'_1)$.

Now assume the statement of the lemma holds for some 
$m \ge 1$, and we prove it for $m+1$. Hence assume that 
(A)--(D) hold with $m$ replaced by $m+1$.
By the induction hypothesis, 
$\eta_m := C'_1 \vee \dots \vee C'_m \in \OmegaL_{1,m}$.
Since the definitions of $B_k$, $f_k$ ($1 \le k \le m$)
do not depend on $C'_{m+1}$, we also get $B_k = B'_k$ and $f_k = f'_k$
for $1 \le k \le m$. Note also that $T_k$, $1 \le k \le m$,
has the same value whether we consider the burning of $\eta_m$
or $\eta_{m+1}$. Since $B_m = B'_m$, we have
\eqnst
{ B_{m+1} \supset g(B_m, C'_{m+1}, \es) 
  = g(B'_m, C'_{m+1}, \es) \not= \es, }
by (B). Therefore, $B_{m+1}$ is not empty. We show that 
this implies that $\eta_{m+1}$ is left-burnable.

First, note that rung $m + 2$ (the ghost rung) can
be burnt. Our assumption (C) says that $f'_{m+1} (G)$
is determined via the function $\phi$ by the data: 
$f'_m = f_m$, $B'_m = B_m$ and $C'_{m+1} = C_{m+1}$, and 
that its value can be obtained as the result of the 
computation in \eqref{e:fk+1}. After the ghost rung has 
been burnt, the burning of rungs $m$ and $m+1$ will follow
the pattern of \eqref{e:fk+1}, with $A := G$. Since the 
computation will stabilize with result 
$f'_{m+1} (G) = G$, this means that eventually everything in
rung $m + 1$ burns. By left-burnability of $\eta_m$,
this means that also all the rungs $1, \dots, m$ burn,
and hence $\eta_{m+1}$ is left-burnable.

By Lemma \ref{lem:whatpairs} (c) and (C), we now have
\eqnst
{ (B_{m+1},f_{m+1}) 
  = \phi(B_m,C_{m+1},f_m)
  = \phi(B'_m,C'_{m+1},f'_m) 
  = (B'_{m+1},f'_{m+1}). }
This advances the induction, and the lemma follows.
\end{proof}

\emph{Proof of Theorem \ref{thm:coding}.}
Let $\psi$ and $\phi$ be as in Lemma \ref{lem:whatpairs}.
We define the inclusion $i(C) = \big(C, \psi(C)\big)$.
We define the alphabet $\cA$ as the set of 
$(C,B,f)$ such that there exists $m \ge 1$ and a sequence
$(C_k,B_k,f_k)_{k=1}^m$ with 
\begin{itemize}
\item[(i)] $(C_m,B_m,f_m) = (C,B,f)$;
\item[(ii)] $(C_1,B_1,f_1) = i(C_1)$; 
\item[(iii)] $g(B_k,C_{k+1},\es) \not= \es$, $1 \le k < m$; 
\item[(iv)] $(B_{k+1},f_{k+1}) = \phi(B_k,C_{k+1},f_k)$, 
  $1 \le k < m$; and 
\item[(v)] $f_k(G) = G$, $1 \le k \le m$. 
\end{itemize}
We define the transition matrix $T$ by
\eqnst
{ T\big((C,B,f),(C',B',f')\big)
  = \begin{cases}
    1 & \text{if } g(B,C',\es) \not= \es,\, (B',f') = \phi(B,C',f)\,, \\
    0 & \text{otherwise.}
    \end{cases} }
It follows from these definitions that for any left-burnable 
$\eta$, 
\eqnst
{ \omega 
  := (\omega_k)_{k=1}^m 
  := (C_k,B_k,f_k)_{k=1}^m \in \cA_{1,m}. }
By Lemma \ref{lem:whatpairs} we have in fact 
$\omega \in \bOmega^{\rm L}_{1,m}$. It follows from 
Lemma \ref{lem:burnable} that every element of 
$\bOmega^{\rm L}_{1,m}$ arises this way. The correspondence
$\eta \mapsto \omega$ satisfies 
$P(\omega) = \eta$ by defintion, and hence is one-to-one.

It remains to show that $T$ is transitive. For this, we first
show that for any $(C,B,f) \in \cA$ we have 
$T\big((C,B,f),i(\Cmax)\big) = 1$. It is easy to verify that 
$i(\Cmax) = (\Cmax,G,\fmax)$, where $\fmax \equiv G$.
By the definition of $\cA$, $B \not= \es$, and hence
$g(B,\Cmax,\es) = G \not= \es$. Recalling the construction 
of $\phi$ in \eqref{e:Bk+1}--\eqref{e:fk+1}, we have, 
regardless of the values of $B$ and $f$, 
\eqnst
{ A^{(0)} = g(B,\Cmax,\es) = G, }
and therefore $A^{(j_0)} = G$. This implies that 
$\bar{B}^{(0)} = f(G) = G$, and hence 
$\bar{A}^{(0)} = g(G, \Cmax, A) = G$, regardless
of what $A$ is. It follows that 
\eqnst
{ \phi(B,\Cmax,f) = (G,\fmax), } 
as required.

Next we show that $T\big(i(\Cmax),i(C)\big) = 1$ for every $C \in \cC$.
Since $i(\Cmax) = (\Cmax,G,\fmax)$, the requirement that 
$g(G,C,\es) \not= \es$ is clearly satisfied. Recalling the 
construction of $\phi$ in \eqref{e:Bk+1}--\eqref{e:fk+1}, 
we have for any $A \subset G$,
\eqnsplst
{ A^{(0)} 
  &= g(G,C,\es) 
  = A^{(j_0)}\,, \\
  \bar{B}^{(0)} 
  &= \fmax(A^{(j_0)}) 
  = G\,, \\
  \bar{A}^{(1)} 
  &= g(G,C,A) 
  = \bar{A}^{(\bar{j}_0)}. }
This shows that $\phi(G,C,\fmax) = \big(g(G,C,\es), g(G,C,\cdot)\big)$,
as required.

We have shown that $i(\Cmax)$ can be reached from any state,
and any $i(C)$ can be reached from $i(\Cmax)$. By the definition 
of $\cA$, any state can be reached from some $i(C)$, and hence
from $i(\Cmax)$. Using again that $i(\Cmax)$ can follow any state,
we see that no periodicity issue can arise, and hence $T$ is
transitive.
\qed

\medbreak

{\bf Acknowledgements.} We thank Omer Angel for useful conversations.

\end{document}